\documentclass{amsart}

\usepackage{amsfonts}
\usepackage{amsmath,amscd,amssymb,amsthm}
\usepackage{mathrsfs}
\usepackage{enumerate,multicol}
\usepackage{url}
\usepackage{pstricks}
\usepackage{pst-plot}
\usepackage{graphicx}

\theoremstyle{plain}
\newtheorem{theorem}{Theorem}

\newtheorem{lemma}[theorem]{Lemma}

\theoremstyle{definition}
\newtheorem{definition}[theorem]{Definition}
\newtheorem{example}[theorem]{Example}

\theoremstyle{remark}

\newcommand{\x}{ {\bf{x}} }

\newcommand{\Z}{\mathbb{Z}}

\newcommand{\N}{\mathbb{N}}

\begin{document}

\title{Standard Young Tableaux and Lattice Paths}

\author{Shaun V. Ault}
\address{
  Department of Mathematics,
  Valdosta State University,
  Valdosta, GA 31698}
\email{svault@valdosta.edu}

\begin{abstract}
  Using lattice path counting arguments, we reproduce a well known
  formula for the number of standard Young tableaux.  We also produce
  an interesting new formula for tableaux of height $\leq 3$ using 
  the Fourier methods of Ault and Kicey.
\end{abstract}

\keywords{Young Tableaux, Lattice Paths, Generating Functions}

%\subjclass[2010]{}

\maketitle

%%%%%%%%%%%%%%%%%%%%%%%%%%%%%%%%%%%%%%%%%%%%%%%%%%%%%%%%%%%%%%%%%%%%%%%%%%%%%%%%

\section{Definitions and Statement of Theorem}

In this short paper, we use lattice path counting arguments to develop
a generating function for the number of standard Young tableaux of
shape $\lambda$, which we denote $f^{\lambda}$.  From
there, a well-known formula for $f^{\lambda}$ can be derived.  The formula
is not new, and neither do we claim to have discovered a new generating
function.  What is new (as far as can be determined) is the connection between
lattice paths and $f^{\lambda}$.  Moreover, the Fourier
methods of Kicey and the author~\cite{AultKicey} provide an interesting
(albeit not necessarily useful) formula for $f^{\lambda}$ where the
height of $\lambda$ is no greater than 3.

We assume that the reader is familiar enough with the basic ideas
of Young tableaux.  No deep representation theory is required.
The following definitions are fairly standard in the literature
(e.g.~\cite{FH,M98,Stanley-AC}). Fix $r \in \N$.  
A partition $\lambda = (\lambda_1, \lambda_2, \ldots, \lambda_r)$
will have the property that $\lambda_i \geq \lambda_j$ as long as $i \leq j$.
For convenience, we allow $\lambda_i = 0$.  Necessarily, all 
zero parts will occur at the end of the sequence.  If $\lambda_r > 0$, then
we say that $\lambda$ has height $r$. 
Denote the
size of $\lambda$ by $|\lambda| = \lambda_1 + \lambda_2 + \cdots + \lambda_r$.
The number of standard
tableaux of shape $\lambda$ will be denoted by $f^{\lambda}$.
Denote by $\Sigma_r$ the symmetric group on the $r$ letters,
which acts on $r$-tuples $\mathbf{x} = (x_1, x_2, \ldots, x_r)$ by
permuting entries:
\[
  \sigma\mathbf{x} = (x_{\sigma^{-1}(1)}, x_{\sigma^{-1}(2)}, 
  \ldots, x_{\sigma^{-1}(r)}).
\]
For any $r \in \N$, let $\mathbf{r} = (0, 1, 2, \ldots, r-1)$,
and denote by $\mathbf{r}^{\ast}$ the reversed $r$-tuple, $(r-1, r-2, \ldots, 1, 0)$.
Let $x_1, x_2, \ldots, x_r$ be formal commuting variables, 
and let $\mathbf{x} = (x_1, x_2, \ldots x_r)$.  Furthermore, if $\mathbf{m}
= (m_1, m_2, \ldots, m_r)$, then let $\mathbf{x}^{\mathbf{m}} = 
x_1^{m_1}x_2^{m_2} \cdots x_r^{m_r}$.  Let
$V_r$ be the Vandermonde polynomial,
\begin{equation}\label{eqn.Vr}
  V_r = \sum_{\sigma \in \Sigma_r} \mathrm{sgn}(\sigma) \mathbf{x}^{\sigma(\mathbf{r}^*)}
  = \prod_{i < j} (x_i - x_j)
\end{equation}
Finally, let $t_r = \sum_{i=1}^{r} x_i$.
We will be working with Laurent polynomials in the variables
$\{x_1, x_2, \ldots x_r\}$.  For such
a function $f$, let $[f]_{\mathbf{m}}$ be the coefficient of the
term $\mathbf{x}^{\mathbf{m}}$.
The following theorem is already present (at least implicitly) in the 
literature.  For
example, Fulton and Harris develop the tools in \S 4 of~\cite{FH}.
\begin{theorem}\label{thm.OGF}
  The number of standard tableaux of shape $\lambda$ whose height is
  no greater than $r$ is 
  equal to the $\mathbf{x}^{\mu}$ coefficient of $t_r^n V_r$,
  where $n = |\lambda|$ and
  $\mu = \lambda+\mathbf{r}^*$.  Equivalently,
  \begin{equation}\label{eqn.f^lambda}
    f^{\lambda} = \left[  \frac{t_r^nV_r}{\mathbf{x}^{\mathbf{\mathbf{r}^*}}}
    \right]_{\lambda}
  \end{equation}
\end{theorem}

Theorem~\ref{thm.OGF} implies that 
$F_{n,r} = \frac{t_r^nV_r}{\mathbf{x}^{\mathbf{\mathbf{r}^*}}}$ 
is a generating function for the numbers $f^{\lambda}$ where
$\lambda$ is a partition of $n$ having up to $r$ components.
Theorem~\ref{thm.OGF} directly implies the following formula, which is
well-known in the literature~\cite{FH}.
  \begin{equation}
    f^{\lambda} = \left(
        \begin{array}{ccccc}
      && n && \\
        \mu_1 & \mu_2 & \mu_3 & \cdots & \mu_r
    \end{array}
        \right) \prod_{i<j} (\mu_i - \mu_j),
  \end{equation}
where $\mu_k = \lambda_k + r - k$.  

\section{Proof of Theorem~\ref{thm.OGF} by Way of Lattice Paths}

\begin{figure}
  \caption{A portion of the underlying graph
   $\Lambda^3$ showing individual
    standard tableaux. Different arrowheads signify the addition of
    different unit vectors.}\label{fig.YoungGraph}
  \begin{center}
    \scalebox{0.4}{\includegraphics{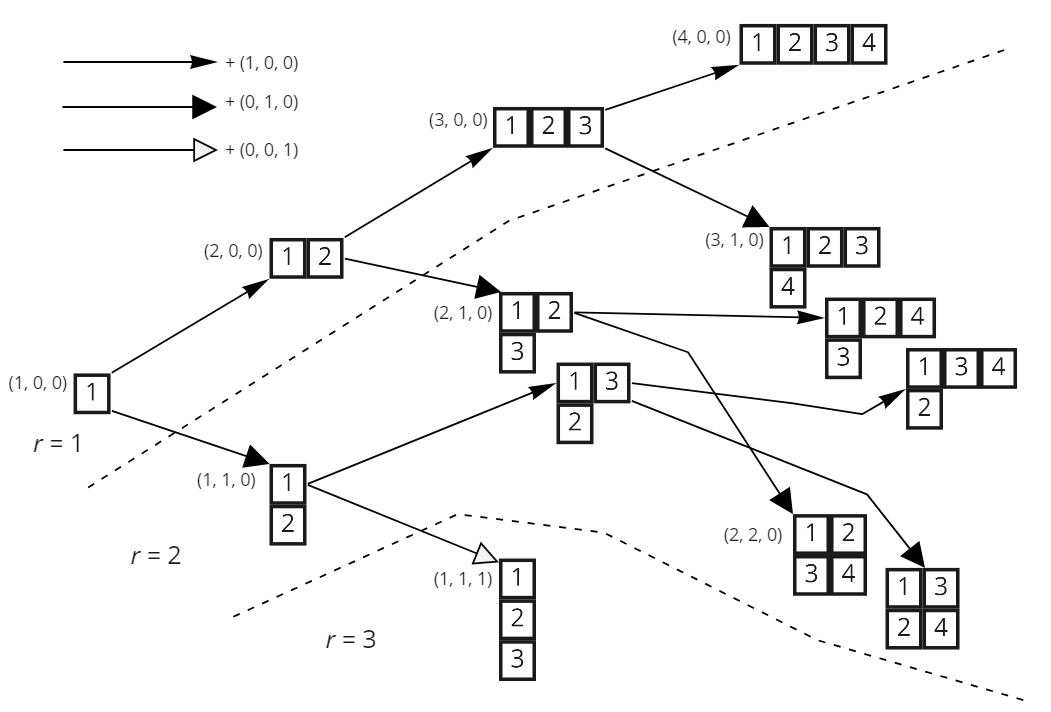}}
  \end{center}
\end{figure}

Fix $r \in \N$ and consider the set,
\[
  \Lambda^r = \{\lambda = (\lambda_1, \lambda_2, \ldots, \lambda_r)
  \in \Z^{r} \;|\;
  \lambda_i \geq \lambda_{j} \geq 0 \;\textrm{for all}\; i \leq j \}.
\]
We may consider elements $\lambda \in \Lambda^r \subseteq \N_0^r$ to 
be partitions, points, or
vectors as necessary. Each $\lambda \in \Lambda^r$
fits into a directed graph determined by the order in which the number
labels are inserted into the Ferrers diagram of the Young tableau.  
For example, for $\lambda = (3, 1)$
there are $f^{\lambda} = 3$ ways to fill in the Young tableau, as suggested
by Figure~\ref{fig.YoungGraph}. (Note, the shape $(3, 1)$ is equivalent to
$(3, 1, 0)$).  Each step in the directed graph adds 1 to
one of the components of the vector, as long as the vector components remain
non-increasing after the addition.  Thus, if $\lambda = (1, 1)$, then
$\lambda + (1, 0) = (2, 1)$ is legal, but $\lambda + (0, 1) = (1, 2)$ is not.

Let
\[
  \delta_i = (0, \ldots, \underset{i}{1}, 
	           \ldots, 0).
\]
There is a recursive formula,
\begin{equation}\label{eqn.f-lambda-recursive}
  f^{\lambda} = \sum_{i = 1}^{r} f^{\lambda - \delta_i},
\end{equation}
where $f^{\mu} = 0$ whenever $\mu$ is not a legal partition -- that is, 
whenever
$\mu \notin \Lambda^r$.  We take the convention that there is exactly
one trivial Young tableau, corresponding to $\lambda = \mathbf{0} =
(0, 0, \ldots, 0)$, so we have $f^{\mathbf{0}} = 1$.

\begin{figure}
  	\caption{$\Lambda^2$ with $\lambda_i \leq 4$, left; 
  	$\Lambda^3$ with $\lambda_i \leq 3$, right.  The vertex numbers count
  	number of paths from the origin to each point.  These numbers
  	coincide with $f^{\lambda}$.}
  	\label{fig.L2-3}
	  	\scalebox{0.34}{\includegraphics{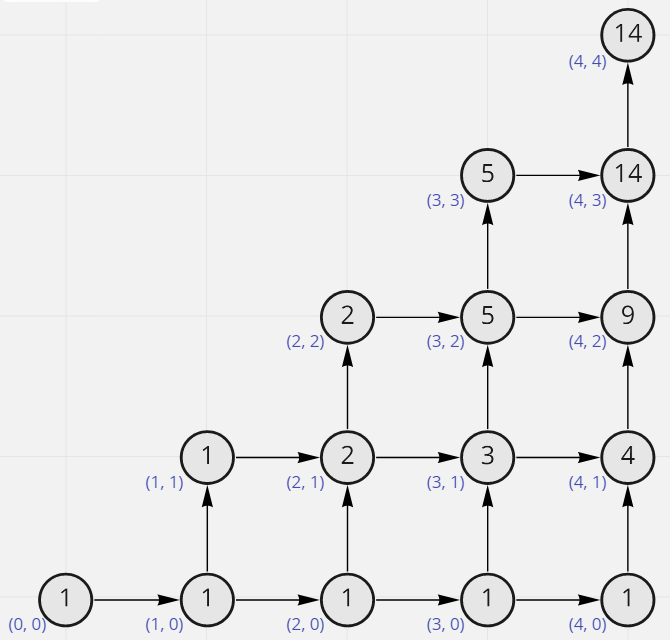}}
	  	\scalebox{0.34}{\includegraphics{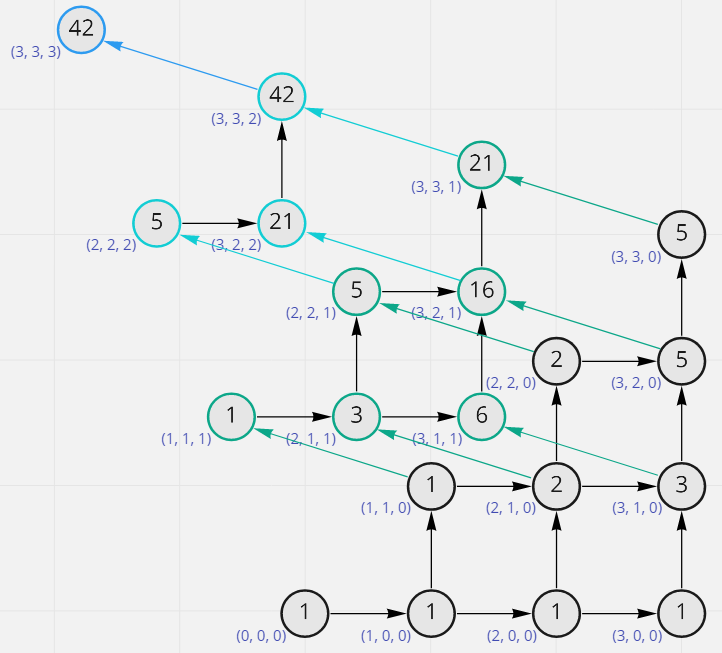}}
\end{figure}

The key to our argument is to interpret partitions
$\lambda$ as points in the lattice $\N_0^r$ with move
set $\{ \delta_i \;|\; 1 \leq i \leq r\}$, where
a move in the direction $\delta_i$ is viewed as 
appending and labeling a new square to the tableau at row $i$.
Thus each path that remains entirely within
$\Lambda^r$ (the legal partitions) beginning at the origin
and ending at the point $\lambda$ represents a distinct way to fill
in the Ferrers digram of shape $\lambda$, and conversely every
standard Young tableau has an associated lattice path within
$\Lambda^r$ determined by its numberic labels.  Therefore, the
count of all such
lattice paths from $\mathbf{0}$ to $\lambda$ is equal to
$f^{\lambda}$.
See Figure~\ref{fig.L2-3} for examples of the lattice $\Lambda^r$ for
$r = 2$ and $r=3$.
  
Counting paths in $\Lambda^r$ is complicated by the restrictions on
legal partitions $\lambda$.  A simple shift of the lattice 
by $\mathbf{r}^*$ suffices
to make things easier.  Indeed, this shift will also allow
Eqn.~(\ref{eqn.f-lambda-recursive})
to be defined more explicitly
in terms of identifying exactly when $\lambda - \delta_i$ is a legal
partition.
Let
\begin{equation}\label{eqn.Lambda^r}
  \widehat{\Lambda}^r = 
  \Lambda^r + \mathbf{r}^{\ast} =
  \{\lambda = (\lambda_1, \lambda_2, \ldots, \lambda_r)
  \in \Z^{r} \;|\;
  \lambda_i > \lambda_{j} \;\textrm{for all}\; i \leq j,\; \textrm{and} \;
  \lambda_{r} \geq 0 \}.
\end{equation}
Observe, if $\lambda \in \Lambda^r$ and $\lambda - \delta_i \notin
\Lambda^r$, then it must necessarily be the case that 
$\lambda - \delta_i + \mathbf{r}*$ has either an entry equal to $-1$ or
a repeated entry.

Inspired by methods found in Ault and Kicey~\cite{AK1,AultKicey},
consider functions $v : \widetilde{\Lambda}^r \to \Z$, 
which may be extended to $v : \N_0^r \to \Z$ 
in a manner to be explained presently.
   
\begin{definition}
   A function $v : \N_0^r \to \Z$ will be called $\Sigma$-admissible
   if it satisfies $v(\sigma \x) = \mathrm{sgn}(\sigma) v(\x)$
   for all $\sigma \in \Sigma_r$.
\end{definition}
   
\begin{lemma}\label{lem.sigma_admissible}
  If $v$ is $\Sigma$-admissible, and if $\x$ has a repeated entry,
  then $v(\x) = 0$.
\end{lemma}

\begin{proof}
   Suppose $x_i = x_j$.  Let $\sigma \in \Sigma_r$ be the transposition
     $(i, j)$.  Then $v(\x) = v(\sigma \x) = -v(\x)$, which implies 
     that $v(\x) = 0$.
\end{proof}

For convenience, we assume all $\Sigma$-admissible functions $v$ take
the value $v(\x) = 0$ if any component of $\x$ is negative.
Thus, if $v : \widehat{\Lambda}^r \to \Z$, then there is a unique, well-defined
extension of $v$ to a function (also called $v$) such that $v : \Z^r \to \Z$
and $v$ is $\Sigma$-admissble; namely, set
$v(\x) = 0$ whenever $\x$ has negative or repeated entries and $v(\sigma \x) = 
\mathrm{sgn}(\sigma) v(\x)$ whenever $\x \notin \widehat{\Lambda}^r$ has
distinct non-negative entries.

\begin{example}
   If we defined $v_0(2, 1, 0) = 1$ (where $r = 3$), and 
   $v_0(\x) = 0$ for all $\x \in \widehat{\Lambda}^3$ for which $\x \neq (2, 1, 0)$,
   then $v_0$ may be extended to be $\Sigma$-admissible by setting
   $v_0(0, 2, 1) = v_0(1, 0, 2) = v_0(2, 1, 0) = 1$,
   $v_0(0, 1, 2) = v_0(2, 0, 1) = v_0(1, 2, 0) = -1$,
   and $v_0(\x) = 0$ on all other points of $\Z^3$.
\end{example}

For each $1 \leq i \leq r$, let $R_i$ ``right-shift'' operator defined 
on functions $v : \Z^r \to \Z$ by
\[
   R_i[v](x_1, \ldots, x_r) = 
   v(\x - \delta_i) =
   v(x_1, \ldots, x_{i}-1, \ldots, x_r).
\]
Let $T_r = \sum_{i = 1}^r R_i$.

\begin{lemma}\label{lem.v-adm}
   If $v$ is $\Sigma$-admissible, then so is $T_r[v]$.
\end{lemma}

\begin{proof}
  It suffices to show that 
  $T_r[v](\sigma\x) = -T_r[v](\x)$
  for a transposition $\sigma$.
  Let $\sigma = (j, k)$ for some $1 \leq j < k \leq r$,
  and let $\x_{jk} = \sigma\x$.  Clearly, $\x_{jk}$ has the
  same entries as $\x$ but with entries $j$ and $k$ swapped.
  \[
    T_r[v](\sigma\x) =
    \sum_{i=1}^{r} R_i[v](\x_{jk}) 
    = \sum_{i=1}^{r}v\left(\x_{jk}  - \delta_i\right)
  \]
  Note that if $i \neq j, k$, then we have by $\Sigma$-admissibility
  of $v$,
  \[
    v(\x_{jk} - \delta_i) = -v(\x - \delta_i)
  \]
  Moreover, for $i = j$ and $i = k$, we have:
  \begin{eqnarray*}
    v(\x_{jk} - \delta_j) &=& -v(\x - \delta_k) \\
    v(\x_{jk} - \delta_k) &=& -v(\x - \delta_j)
  \end{eqnarray*}
  Therefore, $T_r[v](\sigma\x) = \sum_{i=1}^r \left(-v(\x - \delta_i)\right) =
  -T_r[v](\x)$, as required.
\end{proof}

Define $v_0 : \widetilde{\Lambda}^r \to \Z$ such that $v_0(\mathbf{r}^*) = 1$
and $v_0(\x) = 0$ for $\x \neq \mathbf{r}^*$, and extend $v_0$ to
a $\Sigma$-admissible function.
For each $n \in \N$, let $v_n = T_r^n[v_0]$, which by Lemma~\ref{lem.v-adm} 
is $\Sigma$-admissible.
	
\begin{lemma}\label{lem.paths}
  $v_n(\x)$ counts the number of paths of length $n$ from $\mathbf{r}^*$
  to $\x$ that remain entirely in the lattice $\widetilde{\Lambda}^r$,
  using the move set $\{ \delta_i \;|\; 1 \leq i \leq r\}$.
\end{lemma}

\begin{proof}
  There is one $0$-length path beginning and ending at $\mathbf{r}^*$;
  therefore $v_0$ counts all $0$-length paths.
  now fix $n \in \N$ and suppose that $v_n(\x)$ counts all $n$-length
  paths from $\mathbf{r}^*$ to $\x$ within $\widetilde{\Lambda}^r$.
  Because $v_n$ is $\Sigma$-admissible, we have $v_n(\x) = 0$ on all boundary
  points $\x$ adjacent to the lattice.  Thus, for any 
  $\x \in \widetilde{\Lambda}^{r}$, we have
  \[
    v_{n+1}(\x) = T[v_n](\x) = \sum_{i=1}^r v_n(\x - \delta_i),
  \]
  where $v_n(\x - \delta_i) = 0$ for any case in which $\x - \delta_i
  \notin \widetilde{\Lambda}^r$.  Thus $v_{n+1}$ counts all paths
  from $\mathbf{r}^*$ to $\x$ of length $n+1$ that remain entirely within
  $\widetilde{\Lambda}^r$.
\end{proof}	

\begin{lemma}\label{lem.lambda_v}
   If $|\lambda| = n$, then
  	$v_n(\lambda + \mathbf{r}^*) = f^{\lambda}$.
\end{lemma}
   
\begin{proof}
   Let $\mathbf{0} = (0, 0, \ldots, 0)$ be the empty partition.
   By~(\ref{eqn.Lambda^r}) and Lemma~\ref{lem.paths}, we have
   $v_0(\mathbf{0} + \mathbf{r}^*) = 1 = f^{\mathbf{0}}$.
   Next, assume that
   Lemma~\ref{lem.lambda_v} is true for all $\lambda$ of length
   $|\lambda| = n$ (for some fixed $n \geq 0$).
   Making use of Eqn.~(\ref{eqn.f-lambda-recursive}), we may
   derive the formula for $\lambda$ with $|\lambda| = n+1$ as follows.
   Lemma~\ref{lem.sigma_admissible} shows that
   $v_k$ is zero on the boundaries where 
   $\lambda - \delta_i + \mathbf{r}^* \notin
   \widetilde{\Lambda}^r$, which is equivalent to $\lambda - \delta_i
   \notin \Lambda^r$; in other words, whenever $\lambda - \delta_i$
   is not a legal partition due to having negative or repeated entries. 
   Then it follows that:    
   \[
     f^{\lambda} = \sum_{i = 1}^{r} f^{\lambda - \delta_i} 
     = \sum_{i = 1}^{r} v_n(\lambda - \delta_i + \mathbf{r}^{*})
     = T_r[v_n](\lambda + \mathbf{r}^*) = v_{n+1}(\lambda + \mathbf{r}^*).
   \]
\end{proof}
   
Finally, we make the connection to Eqn.~(\ref{eqn.f^lambda}).
Recall $V_r$ from Eqn.~(\ref{eqn.Vr}).  
There are exactly $r!$ nonzero terms, including
$\x^{\mathrm{r}^*}$, whose coefficient is 1.
It should be clear from its form as a sum of terms of the form
$\mathrm{sgn}(\sigma)\x^{\sigma(\mathbf{r^*})}$ that
$V_r$ is the (ordinary) generating function corresponding to $v_0$.
Moreover, the transition operator $T_r$ corresponds to multiplication
by $t_r = x_1 + \cdots + x_r$.  Thus, the generating function
that counts the number of $n$-length
paths in $\widetilde{\Lambda}^r$ from
$\mathbf{r}^*$ to $\x$ is precisely $t_r^nV_r$.  That is, $t_r^nV_r$ is
the generating function for $v_n$.
Lemma~\ref{lem.lambda_v} then implies that a shift of exponents
is needed to get the generating function for $f^{\lambda}$; namely, 
reduction of all multi-exponents by $\mathbf{r}^*$, which corresponds
to division by $\x^{\mathbf{r}^*}$ in the generating function.  This completes
the argument and proves Theorem~\ref{thm.OGF}.

\section{Examples}

\begin{example}
  Counting $f^{(k, \ell)}$, or standard tableaux with at most 
  two rows and $k \geq \ell$.
  Here, $\mathbf{r} = (0, 1)$, $V_r = V_2 = x_1 - x_2$, and $t_r = t_2 = x_1 + x_2$.
  \[
    F_{n,2} = \dfrac{t_2^n V_2}{x_1^1 x_2^0} = \dfrac{(x_1+x_2)^n(x_1 - x_2)}{x_1}
    = (x_1 + x_2)^n\left(1- \frac{x_2}{x_1}\right)
  \]
  The first few of these are:
  \begin{eqnarray*}
     F_{1,2} &=& x_1 - \frac{x_2}{x_1} \\
     F_{2,2} &=& x_1^2 + x_1x_2 - \frac{x_2^3}{x_1} \\
     F_{3,2} &=& x_1^3 + 2x_1^2x_2 - 2x_2^3 - \frac{x_2^4}{x_1} \\
     F_{4,2} &=& x_1^4 + 3x_1^3x_2 + 2x_1^2x_2^2 - 2x_1x_2^3 - 3x_2^4 - 
     \frac{x_2^5}{x_1} \\
     F_{5,2} &=& x_1^5 + 4x_1^4x_2 + 5x_1^3x_2^2 - 5x_1x_2^4 - 4x_2^5 - 
     \frac{x_2^6}{x_1}\\
     F_{6,2} &=&
     x_1^6 + 5x_1^5x_2 + 9x_1^4x_2^2 + 5x_1^3x_2^3 - 5x_1^2x_2^4 - 9x_1x_2^5
      - 5x_2^6 - \frac{x_2^7}{x_1}\\
     F_{7,2} &=&     
     x_1^7 + 6x_1^6x_2 + 14x_1^5x_2^2 + 14x_1^4x_2^3 - 14x_1^2x_2^5 - 14x_1x_2^6
      - 6x_2^7 - \frac{x_2^8}{x_1}
  \end{eqnarray*}
  
  Note that the terms whose exponents are multi-indeces that correspond to 
  legal partitions $\lambda$ have the expected coefficients.  For example,
  we can look at $F_{7,2}$.
  \[
    [ F_{7,2} ]_{(7,0)} = 1 = f^{(7,0)} \quad
    [ F_{7,2} ]_{(6,1)} = 6 = f^{(6,1)} 
    \]
    \[
    [ F_{7,2} ]_{(5,2)} = 14 = f^{(5,2)} \quad
    [ F_{7,2} ]_{(4,3)} = 14 = f^{(4,3)} 
  \]
  
  In this simple case, the Binomial Theorem can be used to
  produce an explicit formula.
  \begin{eqnarray*}
    F_{n, 2} &=& \left(1 - \frac{x_2}{x_1}\right)
    \sum_{k=0}^n \binom{n}{k} x_1^k x_2^{n-k}\\
    &=& \sum_{k=0}^n \binom{n}{k} x_1^k x_2^{n-k}
    - \sum_{k=0}^n \binom{n}{k} x_1^{k-1} x_2^{n-k+1} \\
    &=& \sum_{k=0}^n \binom{n}{k} x_1^k x_2^{n-k}
    - \sum_{k=-1}^{n-1} \binom{n}{k+1} x_1^{k} x_2^{n-k} \\
  \end{eqnarray*}
  Thus, with $\ell = n - k$, if $k \geq \ell$, then the above yields:
  \[
    f^{(k, \ell)} = \binom{k+\ell}{k} - \binom{k+\ell}{k+1}
  \]  
\end{example}

\section{Slices of the Lattice}

The functions $v_n(\x)$ defined above may be interpreted in a 
different way -- as walks in Weyl alcoves.  
For $n \in \N_0$, let $\Lambda^r_n = \{\lambda \in \Lambda \;|\; |\lambda| = n\}$,
which we call the $n$th slice of $\Lambda^r$.
Observe that $\Lambda^r_n$ is a subset of the
plane $x_1 + x_2 + \cdots + x_r = n$.  This plane can be viewed as an affine
lattice $A_{n-1}$,
whose move set is determined by the simple roots $\alpha_i = \delta_i - \delta_{i+1}$.
Then $\Lambda_n^r$, with an appropriate shift,
exists within the interior points of an alcove in $A_{n-1}$.
 
\begin{figure}
\caption{Lattice $A_2$.  Three vectors are highlighted, $\alpha$,
$\beta$, and $-\alpha - \beta$, which we take as the permissible
move set.}\label{fig.A2}
\begin{center}
  \scalebox{0.75}{\includegraphics{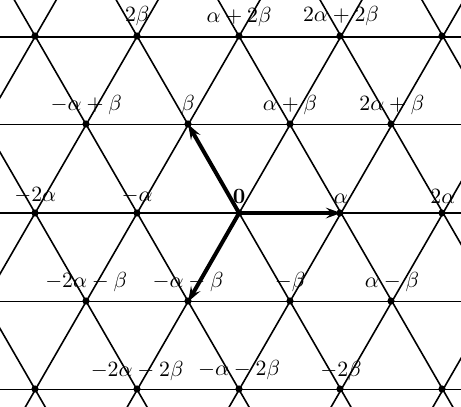}}
\end{center}
\end{figure}

We shall focus on the case $r = 3$, 
that is, we shall develop a formula to count $f^{\lambda}$ where
$\lambda = (a, b, c)$,
as this is precisely the
case for which methods of Ault-Kicey~\cite{AultKicey} will be
useful and interesting.  (Unfortunately, the Fourier methods of
Ault-Kicey do not easily scale beyond $A_2$ lattices.  For more about
walks in arbitrary Weyl alcoves and related structures, 
see~\cite{GesselKrattenthaler,GesselZeilberger,Grabiner,MortimerPrellberg}.)
Let $A_2$ be defined in the usual way as a lattice in the plane
generated by roots $\alpha = (1, 0)$ and $\beta = (-1/2, \sqrt{3}/2)$,
as shown in Figure~\ref{fig.A2}.
Thus, all points of $A_2$ take the form $\langle u, v \rangle =
u\alpha + v\beta$ for $u, v \in \Z$.

Consider the
affine transformation,
\begin{eqnarray}
  P &:& \Lambda_n^3 \to A_2 \nonumber \\
  P(x, y, z) &=& \langle x-z+2, y-z+1 \rangle \label{eqn.P}\\
             &=& \begin{bmatrix}
                    1 & 0 & -1 \\ 0  & 1 & -1
                    \end{bmatrix}
                    \begin{bmatrix} x \\ y \\ z \end{bmatrix}
                    + \begin{bmatrix} 2 \\ 1 \end{bmatrix}
\end{eqnarray}
With $u = x-z + 2$ and $v = y-z+1$ as defined by Eqn.~(\ref{eqn.P}),
we have $u > v$ and $v > 0$.  This implies that $P(x, y, z)$ lies
in the interior of the region bounded by the $\alpha$-axis and
$\alpha+\beta$-axis.  We also find that since $x \leq n$, then
$u = x - z + 2 \leq n + 2$.  That is, $P(\Lambda_n^3)$ lies
entirely within the triangular region determined by $u > v$,
$v > 0$, and $u < n+3$ in $A_2$.  

Moreover, the three moves, $\delta_1 = (1, 0, 0)$, $\delta_2 = (0, 1, 0)$,
and $\delta_3 = (0, 0, 1)$, in $\Lambda^3$ transform to the
following moves in $A_2$:
\[
  \begin{bmatrix}
     1 & 0 & -1 \\ 0  & 1 & -1
  \end{bmatrix}\delta_1 = \langle 1, 0\rangle = \alpha, \quad
  \begin{bmatrix}
     1 & 0 & -1 \\ 0  & 1 & -1
  \end{bmatrix}\delta_2 = \langle 0, 1\rangle = \beta,
\]
\[
  \begin{bmatrix}
     1 & 0 & -1 \\ 0  & 1 & -1
  \end{bmatrix}\delta_3 = \langle -1, -1\rangle = -\alpha - \beta
\]
Thus, counting paths in $\Lambda^3$ is equivalent to counting
paths in a triangular lattice of sufficiently large size, say 
$a + 1$ points along a side 
(where $\lambda = (a, b, c)$ is the shape for which
we wish to count $f^{(a, b, c)}$, where $a + b + c = n$) with
move set $\{ \alpha, \beta, -\alpha-\beta \}$.
Refering to the notation and methods established in
Appendix B of~\cite{AultKicey}, 
the relevant transition operator 
is $T^+ = R^{\alpha} + R^{\beta} + R^{-\alpha-\beta}$, and so we
have on the Fourier side,
\[
  \widehat{T^+}(\omega_1, \omega_2) = e^{-\frac{2\pi i \omega_2}{3(a+1)}}
  + e^{-2\pi i \left(\frac{\omega_1}{a+1} - \frac{2\omega_2}{3(a+1)}\right)}
  + e^{-2\pi i \left(-\frac{\omega_1}{a+1} + \frac{\omega_2}{3(a+1)}\right)}
\]
The Fourier-transformed initial state is defined by
\[
  V_0(\omega_1, \omega_2) = 2i\left[
  -\sin\left(
    \frac{2\pi \omega_1}{a+1}\right)
  - \sin\left(
    2\pi \left(\frac{\omega_1 - \omega_2}{a+1}\right)\right)
  + \sin \left(2\pi \left(\frac{2\omega_1 - \omega_2}{a+1}\right)\right)
  \right]
\]
Now, according to Theorem~B.1 of~\cite{AultKicey},
$v_n = \mathcal{F}^{-1}\left[\widehat{T^+}^nV_0\right]$ 
is the vertex function that counts the paths in the wedge that we
are interested in.  Here $\mathcal{F}$ is a discrete Fourier
transform, and $\mathcal{F}^{-1}$ is its inverse.
Adjusting the input by the affine transformation
Eqn.~(\ref{eqn.P}), we have produced an interesting formula for
height-3 standard Young tableaux.
\begin{equation}\label{eqn.height3}
  f^{(a, b, c)} = \frac{1}{3(a+1)^2}
  \sum_{\omega_1=0}^{a}
  \sum_{\omega_2=0}^{3a+2}
  e^{2\pi i \left( \frac{(a-c+2)\omega_1}{a+1} + \frac{(b-c+1)\omega_2}{3(a+1)}
  \right)}
  \widehat{T^+}(\omega_1, \omega_2)
  V_0(\omega_1, \omega_2)
\end{equation}

The author concedes that Eqn.~(\ref{eqn.height3}) may be of very
little practical use to those who actually wish to count $f^{\lambda}$,
except perhaps as a curiosity.
However, the formula is fairly easy to code into a computer
algebra system such as Sage~\cite{sage}, and has been verified
to produce correct counts of $f^{\lambda}$.

\section{Acknowledgements}

The investigation into partitions and Young tableaux was inpired by
a recent talk by Robert (Bob) Donley~\cite{DonleyTalk}.  
Bob\footnote{Bob Donley teaches as the
Queensborough Community College (New York, NY) and is an organizer of the 
Representation Theory Seminar at the CUNY Graduate Center.} provided
many insightful suggestions that helped to shape the direction of this paper.
I am also grateful to Cyril Banderier and the other organizers
of the virtual conference, ``Lattice Paths, Combinatorics and
Interactions,'' hosted by the Centre International de Rencontres Mathématiques
(June 2021) for accepting my poster presentation, ``From Lattice Paths
to Standard Young Tableaux.''

%%%%%%%%%%%%%%%%%%%%%%%%%%%%%%%%%%%%%%%%%%%%%%%%%%%%%%%%%%%%%%%%%%%%%%%%%%%%%%%%

\bibliographystyle{plain}
\bibliography{refs} %% Use with bibtex

\end{document}